\documentclass[12pt]{amsart}
\usepackage{amscd,graphicx}
\theoremstyle{definition}

\def \ph{\varphi}

\newcommand{\out}{\operatorname{Out}}
\renewcommand{\P}{\mbox{$\mathbb P$}}

\def \ra{\rightarrow}

\def \hom{\mbox{\rm Hom}}
\def \ie{\hbox{\it i.e.}}

\def \sl{\mbox{$\mathfrak{sl}$}}
\def \tns{\otimes}

\def \mwedge{\wedge\cdots\wedge}

\def \mcom{,\cdots,}

\def \C{\mbox{$\mathbb C$}}
\def \R{\mbox{$\mathbb R$}}
\def \Z{\mbox{$\mathbb Z$}}

\def\sh{\operatorname{Sh}}

\def\ad{\operatorname{ad}}

\def\inv{^{-1}}

\def\im{\operatorname{Im}}
\def\A{\mbox{$\mathcal A$}}
\def\B{\mbox{$\mathcal B$}}
\def\L{L}

\def\m{\mbox{$\mathfrak m$}}

\def\tens{\bigotimes}

\def\V{V}

\def\diag{\operatorname{Diag}}
\def\coder{\operatorname{Coder}}

\def\linf{\mbox{$L_\infty$}}
\def\and{\mbox{ \rm and }}

\def\EV{\bigwedge(V)}
\def\aut{\operatorname{Aut}}
\def\dinf{d^{\text{ inf}}}
\def\dinfty{d^{\infty}}
\def\htens{\widehat{\tens}}
\def\m{\mathfrak m}

\def\s#1{(-1)^{#1}}
\DeclareMathOperator*{\invlim}{\overleftarrow{\rm lim}}

\def\psa#1#2{\psi^{#1}_{#2}}

\def\phd#1#2{\ph^{#1}_{#2}}

\def\inv{^{-1}}

\author{Carolyn Otto}
\address{University of Wisconsin\\
Eau Claire, WI 54702-4004} \email{ottoca@uwec.edu}
\author{Michael Penkava}
\address{University of Wisconsin\\
Eau Claire, WI 54702-4004} \email{penkavmr@uwec.edu}
\thanks{The research of the authors was partially supported by
grants from NSF and
 grants from the University of Wisconsin-Eau Claire}
\title[The Moduli Space of 3-d Lie Algebras]
{The Moduli Space of Three Dimensional Lie Algebras}
\begin{document}
\setlength{\multlinegap}{0pt}
%\nocite{ps2,pen1,ls,pen2,pen3,kon,
%fm,mar,mar2,ksv,sta1,sta2,sta3,getz,getz2,ge_ka1,ge_ka2,
%gers,lod,hoch,fi,fi2,ff2,ff,ff3}
%\input abstract.tex
\begin{abstract}
In this paper, we consider the versal deformations of three
dimensional Lie algebras.  We classify Lie algebras and study their
deformations by using linear algebra techniques to study the
cohomology. We will focus on how the deformations fasten the space
of all such structures together. This space is known as the moduli
space. We will give a geometric description of this space, derived
from deformation theory, in order to illustrate general features of
Lie algebras' moduli spaces.
\end{abstract}
\date\today
\maketitle
%\table
%\input intro.tex
%\input part1.tex
\section{Introduction}
Here we establish the basic language in which we will express the
Lie Algebras and their cohomology.  Recall that the exterior algebra
$\bigwedge V$ of an ungraded vector space $V$ has a natural
\Z-grading, with $\deg(v_1\mwedge v_n)=n$ and that there is a
corresponding \Z-graded coalgebra structure with comultiplication
$\Delta$ given by
\begin{equation*}
\Delta(v_1\wedge\cdots\wedge
v_n)=\sum_{k=1}^{n-1}\sum_{\sigma\in\sh(k,n-k)}\s{\sigma}
v_{\sigma(1)}\wedge\cdots\wedge v_{\sigma(k)}\tns
v_{\sigma(k+1)}\wedge\cdots\wedge v_{\sigma(n)}.
\end{equation*}
Here $\sh(k,n-k)$ represents the unshuffles of type $(k,n-k)$, that
is, the permutations which are increasing on $1,\dots,k$ and
$k+1,\dots,n$ and $\s{\sigma}$ represents the sign of the
permutation $\sigma$.   A linear map $\ph\in\hom(\bigwedge^k V,V)$
has degree $k-1$, and extends uniquely to a coderivation by the rule
\begin{equation*}
\ph(v_1\wedge\cdots\wedge
v_n)=\sum_{\sigma\in\sh(k,n-k)}\s{\sigma}\ph(v_{\sigma(1)}\wedge\cdots
v_{\sigma(k)})\wedge v_{\sigma(k+1)}\mwedge v_{\sigma(n)}.
\end{equation*}
Recall that a coderivation of $\bigwedge V$ is a linear map
$\ph:\bigwedge V\ra\bigwedge V$ satisfying
$\Delta\ph=(\ph\tns1+1\tns\ph)\Delta$. We can view
$L_n=\hom(\bigwedge^n V,V)$ as a subspace of the space
$L=\coder(\bigwedge V)$.  Moreover $L=\prod_{n=1}^\infty\L_n$ in a
natural way.  The space of coderivations of $V$ has a natural
structure of a \Z-graded Lie algebra. If $\ph\in\L_k$ and
$\psi\in\L_l$ then $[\ph,\psi]\in\L_{k+l-1}$  is given by $
[\ph,\psi]=\ph\psi-\s{\deg\ph\deg\psi}\psi\ph. $ More explicitly, we
compute
\begin{eqnarray*}
\lefteqn{(\ph\psi)(v_1\mwedge v_{k+l-1})=}&\\
&\sum_{\sigma\in\sh(k,l-1)}\s{\sigma}\ph(\psi(v_{\sigma(1)}\mwedge
v_{\sigma(l)})\wedge v_{\sigma(l+1)}\mwedge v_{\sigma(k+l-1)}),
\end{eqnarray*}
which determines $\ph\psi$ as an element of $L_{k+l-1}$.

If $V=\langle f_i, i=1\dots N \rangle$, and $I=(i_1\mcom i_k)$ is an
increasing multi-index; \ie, $1\le i_1<\cdots < i_k\le N$, then
$\ph^I_k\in L_k$ is defined by $\phi^I_k(f_J)=\delta^I_J f_k$, where
$f_J=f_{j_1}\mwedge f_{j_k}$.  Then we have
\begin{equation*} L_k=\langle \ph^I_j|1\le i_1\le\cdots\le i_k\le N, 1\le j\le N\rangle,
\qquad \dim(L_k)=\binom nk.
\end{equation*}
When $k$ is odd, the elements in $L_k$ are called even and when $k$
is even the elements in $L_k$ are odd. To emphasize the difference,
we will denote the basis elements of $L_{2k}$ in the form $\psi_j^I$
rather than $\ph_j^I$. It is not difficult to give an explicit
formula for the brackets of the coderivations in terms of the basis
$\ph_j^I$. An element $d$ in $L_2$ is called a codifferential if
$d^2=0$, or equivalently $[d,d]=0$, since $[d,d]=2d^2$. The
condition $d^2=0$ is the Jacobi relation
\begin{equation*}
d(d(a,b),c)-d(d(a,c),b)+d(d(b,c),a)=0,
\end{equation*}
so a codifferential is a Lie algebra structure on $V$.

The coboundary $D:L\ra L$ determined by $d$ is given by
$D(\ph)=[d,\ph]$. The fact that $D$ is a differential on $L$; $\ie,
D^2=0$, follows immediately form the fact that $d$ is a
codifferential. Then $H(d)=\ker(d)/\im(d)$ is called the cohomology
of $d$. Since $[d,L_k]\subseteq L_{k+1}$, we can define the $k^{th}$
cohomology group
\begin{equation*}
H^k(d)= \ker(d:L_k\ra L_{k+1})/\im(d:L_{k-1}\ra L_k).
\end{equation*}
Then $H(d)=\prod_{k=1}^\infty H^k(d)$.

It is useful to have an explicit formula for the $n^{th}$ element is
an ordered basis of $\bigwedge^k V$.  If we define $S(n,k)$
recursively as follows by
\begin{align*}
S(n,1)&=(n)\\
S(n,k)&=S(n-\binom {l-1}k,k-1),l),\qquad \binom {l-1}k < n \le
\binom lk.
\end{align*}
Then $f_{S(n,k)}$ is the $n^{th}$ basis element of $\bigwedge^kV$,
if $1\le n\le\binom Nk$.

Any coderivation in $L_k$ can be expressed in the form
\begin{equation}
\ph=a_j^i\ph^{S(j,k)}_i,
\end{equation}
using the Einstein summation convention, so that we can represent
$\ph$ but the $N\times \binom Nk$ matrix $A=a^i_j$. More generally,
$\ph :L_n \ra L_{m-k+1}$ is responsible by the $\binom N{m-k+1}
\times \binom Nm$ matrix $B=b^i_j$, where
\begin{equation*}
\ph(f_{S(j,m)})=b^i_j f_{S(i,m-k+1)}.
\end{equation*}
In particular when $m=2k-1$ the matrix $AB$ represents $\ph^2$. Thus
$\ph^2=0$ is equivalent to the matrix equation $AB=0$. Moreover, if
$d\in\L_2$ is a codifferential, then the Jacobi relation is
equivalent to $AB=0$ where $d=a_j^i \ph^{S(j,2)}_i$ and the
coefficients $b_j^i$ are determined as follows.

It is easy to see that
\begin{eqnarray*}
\lefteqn{\ph^{(u,v)}_w(f_{(r,s,t)})=}&\\& [\delta^u_r
\delta^v_s(\delta^k_w \delta^l_t-\delta^k_t \delta^l_w) -\delta^u_r
\delta^v_t(\delta^k_w \delta^l_s-\delta^k_s \delta^l_w)+ \delta^u_s
\delta^v_t(\delta^k_w \delta^l_r-\delta^k_r \delta^l_w)]f_{(k,l)}.
\end{eqnarray*}
Suppose that
\begin{align*}
S(j,3)&=(r,s,t)\\
S(i,2)&=(k,l).
\end{align*}
Then
\begin{equation*}
b^i_j=a^w_{S\inv(u,v)}[\delta^u_r \delta^v_s(\delta^k_w
\delta^l_t-\delta^k_t \delta^l_w) -\delta^u_r \delta^v_t(\delta^k_w
\delta^l_s-\delta^k_s \delta^l_w)+ \delta^u_s \delta^v_t(\delta^k_w
\delta^l_r-\delta^k_r \delta^l_w)].
\end{equation*}
The above formula is easily implemented on a computer and determines
the coefficients of $B$ linearly in terms of the coefficients of
$A$. Since $AB$ is an $N \times \binom N3$ matrix, this means the
Jacobi relation is given by $N \binom N3$ homogenous quadratic
equations in the $N \binom N2$ coefficients of $a^i_j$.

An invertible linear map $g:V\ra V$ extends to a coalgebra
automorphism of $\bigwedge V$ by $g(v_1\mwedge v_n)=g(v_1)\mwedge
g(v_n)$. Two Lie algebra structures are isomorphic, or in other
words, they are associated codifferentials $d,d'$ are equivalent if
there is an invertible linear map $g:V\ra V$ such that
$g^*(d)=g^{-1}\circ d \circ g =d'$. The set of equivalence classes
of codifferentials on $V$ is called the moduli space of Lie algebra
structures on $V$.

Let $G=g^i_j$ where $g(f_j)=g^i_jf_i$ represents $g:V\ra V$ as an $N
\times N$ matrix. Define the $\binom N2 \times \binom N2$ matrix
$Q=q^i_j$ by $g(f_{S(j,2)})=q^i_jf_{S(i,2)}$ so that $Q$ represents
$g:\bigwedge^2V \ra \bigwedge^2V$. Suppose that
\begin{align*}
S(j,2)&=(u,v)\\
S(i,2)&=(k,l).
\end{align*}
Then
\begin{equation*}
q^i_j=g_u^kg_v^l-g^l_ug^k_v.
\end{equation*}

If $A,A'$ represent the codifferential $d,d'$ respectively, the
condition $g^*(d)=d'$ is represented by the matrix equation $A=G\inv
AQ$. It is easier to solve the equation $GA'=AQ$ as long as the
solution matrix $G$ satisfies $\det G \ne0$.
\section{Three Dimensional Lie Algebras}
Let $V=\langle f_1,f_2,f_3 \rangle$ be a three dimensional vector
space with Lie algebra structure determined by the codifferential
$d$ which is represented by the $3\times 3$ matrix $A=(a^i_j)$. Then
if
\begin{equation*}
B=\left[ \begin {array}{c} -a_{{1,2}}-a_{{2,3}}\\
\noalign{\medskip}a_{{1,1}}-a_{{3,3}}\\\noalign{\medskip}a_{{2,1}}+a_{{3,2}}\end
{array} \right],
\end{equation*}
the Jacobi identity is equivalent to the matrix equation $AB=0$.

Recall that the derived subalgebra of a Lie algebra is the image of
$d:S^2(V) \ra V$. Since $A$ is a matrix representing $d_1$ it
follows that the rank of $A$ is the dimension of the derived
subalgebra.

We first consider the case where the derived subalgebra has
dimension three. Then the matrix $A$ of $d$ is invertible, so it
must be the case that $B=0$. In \cite{fp8} it is shown that any $n$
dimensional Lie algebra structure such that the bracket of any two
elements is a linear combination of those elements has an abelian
ideal of dimension $n-1$.  We shall see later that whenever there is
an ideal of dimension $n-1$ the rank of the associated matrix can
never be larger than $n-1$. Thus by choosing an appropriate basis we
may assume that $d(f_1f_2)=f_3$. Taking into account the fact that
$B=0$ we see that the matrix $A$ has to be of the form
$$A=\begin{bmatrix}0&x&y\\0&z&-x\\1&0&0\end{bmatrix}$$ where
$x^2+yz\ne0$.

It is easy to see that $d$ is equivalent to a codifferential whose
matrix is of the form
$A'=\left[\begin{smallmatrix}0&0&\eta\\0&\beta&0\\1&0&0\end{smallmatrix}\right]$.
In fact, if $z\ne0$, then the linear automorphism $g$ whose matrix
is
$G=\left[\begin{smallmatrix}0&0&-z\\0&1&x\\1&0&0\end{smallmatrix}\right]$
yields $A'=G^{-1}AQ$,where $\beta=x^2+yz$, $\eta=z$; if $y\ne0$ then
$G=\left[\begin{smallmatrix}0&0&-z\\0&1&x\\1&0&0\end{smallmatrix}\right]$
gives $\beta=y$,$\eta=-x^2-yz$; while
 if both $y$ and $z$ vanish then
$G=\left[\begin{smallmatrix}-1/2&-1/2&0\\1&-1&0\\0&0&1\end{smallmatrix}\right]$
gives $\beta=\eta=x$. It follows that any codifferential whose
matrix is invertible is equivalent to one whose matrix is of the
form
$A=\left[\begin{smallmatrix}0&0&\mu\\0&\lambda&0\\1&0&0\end{smallmatrix}\right]$.
The automorphism determined by $G=\diag(r,s,rs)$ preserves the form
of $A$ with $\lambda \ra \lambda'=r^2\lambda$, $\mu \ra
\mu'=s^2\mu$, which is enough to see that there is only one
equivalence class of codifferentials with invertible matrix over
$\C$. Moveover
$G=\left[\begin{smallmatrix}0&-1&0\\1&0&0\\0&0&1\end{smallmatrix}\right]$
preserves the form as well, with $\lambda \ra \lambda'=-\mu$, $\mu
\ra \mu'=-\lambda$, which shows that there are at most two
nonequivalent such codifferentials over $\R$. Of course, these two
codifferentials correspond to the two nonequivalent compact real
forms of the simple complex Lie algebra $\sl(2,\C)$.

Now suppose that the derived algebra of a Lie algebra on a vector
space $V$ of dimension $n$ has dimension smaller than $n$.  Then
there is an ideal $V'$ in $V$ of dimension $n-1$. This gives rise to
the exact sequence
$$ 0\ra V'\ra V\ra\C\ra 0,$$ where $\C$ is the trivial lie
algebra. Let $d'$ be the induced codifferential on $V'$, given by
$d'(u,v)=d(u,v)$ for $u,v\in V'$, with matrix $A'$ in terms of the
basis $\{f_1,\mcom f_{n-1}\}$ of $V'$. Let $\{f_n\}$ be a basis of
the complementary subspace to $V'$ and define $\rho:V' \ra V'$ by
$\rho(u)=d(u,f_n)$. Then $\rho$ is a derivation of $d'$.  If
$R=(r^i_j)$ is the matrix of $\rho$, given by
 $\rho(f_j)=r^i_jf_i$, then the matrix of $A$ is just
$A=\left[\begin{smallmatrix}A'&R\\0&0\end{smallmatrix}\right]$.

If $f_n'=f_n+v$, where $v\in\V'$, then the induced derivation
$\rho'(u)=d(u,f_n')$ determines an equivalent coderivation to $d$.
Moreover,  $\rho'(u)=\rho(u)+\ad_v(u)$, where $\ad_v(u)=d(u,v)$ is
the inner derivation of $V'$ determined by $v$. Thus we can always
replace $R$ by the matrix $R'$ given by $\rho'$, and we see that the
extensions of $V'$ by $\C$ are determined by the space of outer
derivations $\out(d')$ of the induced Lie algebra structure on $V'$.
If we denote the space of inner derivations of $V'$ by $\ad(d')$,
then $\out(d')=H^1(d')/\ad(d')$.

Now let us specialize this to the three dimensional case of Lie
algebras.  We use the fact that there are, up to equivalence, only
two Lie algebra structures on a two dimensional vector space.
\begin{itemize}
\item[\emph{Case 1.}] The nonabelian Lie algebra, given by $d(f_1,f_2)=f_1$.
\item[\emph{Case 2.}] The abelian Lie algebra, given by $d(f_1,f_2)=0$.
\end{itemize}
Let us study the possible forms for the matrix $R$ in case 1.  We
have
\begin{align*}
r^1_1f_1+r^2_1f_2=&\rho(f_1)=\rho(d(f_1f_2))=d(\rho(f_1),f_2)+d(f_1,\rho(f_2))\\
=&d(r^1_1f_1+r^2_1f_2,f_2)+d(f_1,r^1_2f_1+r^2_2f_2)=r^1_1f_1+r^2_2f_1
\end{align*}
It follows that $r^2_1=r^2_2=0$.  Consider the inner derivations
\begin{equation*}
\rho_i(u)=d(u,f_i).
\end{equation*}
It is easy to see that $\rho=-r^1_1\rho_1+r^1_2\rho_2$.  Thus every
derivation is inner, and we may assume that $R=0$.  But then
$\{f_1,f_3\}$ span an ideal on which $d$ acts as the zero matrix.
Thus, the nonabelian case reduces to the abelian one.

Now assume that $d'=0$. Then $d$ is represented by the matrix
$A=\left[\begin{smallmatrix}0&R\\0&0\end{smallmatrix}\right]$, where
$R$ is any $2\times 2$ matrix. However, if $g'$ is any linear
automorphism of $V'$, given by a matrix $G$, and $g$ is the linear
automorphism of $V$ given by the matrix
$R=\left[\begin{smallmatrix}G&0\\0&1\end{smallmatrix}\right]$, then
the $d$ is represented by the matrix
$A=\left[\begin{smallmatrix}0&R'\\0&0\end{smallmatrix}\right]$,
where $R'=(G')\inv RG'$.  Thus similar matrices $R$ determine
equivalent codifferentials. Moreover,  multiplication of $R$ by any
nonzero constant $\lambda$ also determines an equivalent
codifferential, corresponding to the linear automorphism given by
the matrix $G=\diag(1,1,\lambda)$.

As a consequence, we immediately reduce to the following
possibilities.
\begin{itemize}
\item[\emph{$d(\lambda:\mu)$}:] given by the
matrix
$R=\left[\begin{smallmatrix}\lambda&1\\0&\mu\end{smallmatrix}\right]$.
Note that $R$ is similar to the matrix
$\left[\begin{smallmatrix}\mu&1\\0&\lambda\end{smallmatrix}\right]$,
so that $d(\mu:\lambda)\sim d(\lambda:\mu)$.  Moreover, if we
consider the linear automorphism determined by the diagonal matrix
$G=\diag(t,1,t)$, then $R$ is replaced by the matrix
$\left[\begin{smallmatrix}t\lambda&1\\0&t\mu\end{smallmatrix}\right]$,
so that $d(t\lambda:t\mu)\sim d(\lambda:\mu)$.  As a consequence, we
can view $(\lambda:\mu)\in\P^1$. The similarity of matrices with the
same eigenvalues determines an action of $\Sigma_2$ on $\P^1$, so
the set of codifferentials of this type are parameterized by the
orbifold $\P^1/\Sigma_2$.

\item[\emph{$d_2$}:] given by the matrix
$R=\left[\begin{smallmatrix}1&0\\0&1\end{smallmatrix}\right]$.

\item[\emph{$d_1$}:] given by the matrix
$R=\left[\begin{smallmatrix}0&1\\0&0\end{smallmatrix}\right]$.
\end{itemize}

Together with the codifferential $d_3$ representing the Lie algebra
structure $\mathfrak{sl}_2(\C)$ these codifferentials represent a
complete classification of all three dimensional Lie algebras. The
codifferential $d_1$ represents the Heisenberg algebra
$\mathfrak{n}_3(\C)$, $d_2$ represents the solvable Lie algebra
$\mathfrak{r}_{3,1}(\C)$, $d(0:1)$ represents the Lie algebra
$\mathfrak{r}_2(\C)\oplus\C$, $d(1:1)$ represents the solvable Lie
algebra $\mathfrak{r}_{3}(\C)$, while $d(\lambda:\mu)$ represents
the solvable Lie algebra $\mathfrak{r}_{3,\mu/\lambda}(\C)$ when
$\lambda\ne\mu$.  Note that our alignment of algebras into families
differs slightly from the classical alignment.  If we interchange
$d_2$ with $d(1:1)$ then the alignment would agree, but our reasons
for not making this change will become apparent.

\section{Miniversal Deformations of Lie Algebras}
The classical notion of an \emph{infinitesimal deformation} of a Lie
algebra structure on a vector space $V$ given by a codifferential
$d\in L_2$ is a coderivation
\begin{equation*}
d_t=d+t\psi,
\end{equation*}
where $\psi\in L_2$ which satisfies the Jacobi identity
infinitesimally;\ie, $[d_t,d_t]=0\mod t^2$. The condition is
equivalent to $D(\psi)=0$, where $D:L\ra L$ is the coboundary
operator determined by $d$, given by $D(\psi)=[d,\psi]$. The fact
that $D^2=0$, so that $D$ is a differential on $L$, follows from the
Jacobi identity for $d$, and the homology of this differential is
the cohomology $H(d)$ (with coefficients in the adjoint
representation).  Let $\A=\C[t]/(t^2)$. Then $\A$ is a local
algebra, with maximal ideal $\mathfrak m=(t)$. Then we can consider
$d_t$ as being an element of $L\tns\A$, which may be thought of as
the coderivations of $\EV$ with coefficients in $\A$. In this sense
$[d_t,d_t]=0$, in terms of the natural bracket on $L\tns\A$.

One can generalize this construction by allowing $\A$ to be any
local (commutative) algebra over $\C$, and say that $d_{\A}\in
L\tns\A$ is a deformation of $d$ with \emph{local base} $\A$ if
$\epsilon_*(d_{\A})=d$, where $\epsilon:\A\ra\C$ is the
\emph{augmentation} homomorphism, and $\epsilon_*:L\tns\A\ra
L\tns\C=L$ is the induced map.

If $\mathfrak m^2=0$, then we call $\A$ an \emph{infinitesimal
algebra}, and say that $d_{\A}$ is an \emph{infinitesimal
deformation}. This notion of infinitesimal deformation generalizes
the classical notion. An infinitesimal deformation $d_{\A}$ with
base $\A$ is \emph{universal} if whenever $d_{\B}$ is any
infinitesimal deformation with base $\B$, then there is a unique
homomorphism $f:\A\ra\B$ such that $f_*(d_{\A})\sim d_{\B}$, where
$\sim$ means \emph{infinitesimal equivalence}.  Two deformations
$d_{\B}$ and $d'_{\B}$ with the same base $\B$ are said to be
infinitesimally equivalent if there is an \emph{infinitesimal
automorphism} $g_{\B}\in\aut(V)\tns\B$, satisfying
$g_{\B}^*(d_{\B})=d'_{\B}$.  We say that $g_{\B}$ is an
infinitesimal automorphism if $\epsilon_*(g_{\B})=1_V$.

If $\dim(H^2(d))=n$, then the deformation $\dinf$, with base
$\A=\C[t^1,\dots,t^n]/(t^1,\mcom,t^n)^2$, given by
\begin{equation*}
\dinf=d+\delta_i t^i,
\end{equation*}
where $\{\delta_i\}$ is a prebasis of $H^2(d)$ is a universal
infinitesimal deformation.

The classical notion of a formal deformation is given by a power
series
\begin{equation*}
d_t=d+\psi_i t^i,
\end{equation*}
where $\psi_i\in L_2$ and $[d_t,d_t]=0$. We may consider $d_t\in
L\htens\A$, where $\A=\C[[t]]$ is the ring of formal power series in
the variable $t$, and $\htens$ means the formal completion of the
tensor product $L\htens\A=\displaystyle\invlim_{k\ra\infty}L\tens
A/\m^k$, where $\m=(t)$. More generally, we say that $\A$ is a
\emph{formal algebra} if
$\A=\displaystyle\invlim_{k\ra\infty}A/\m^k$, and say that a
codifferential $d_{\A}\in L\htens\A$ is a \emph{formal deformation}
of $d$ if $\epsilon_*(d_{\A})=d$.  There is no universal formal
deformation of $d$, but there is a \emph{versal deformation}.

A formal deformation $d_{\A}$ with base $\A$ is versal if given any
formal deformation $d_{\B}$ with base $\B$, there is a homomorphism
$f:\A\ra\B$ such that $f^*(d_{\A})\sim d_{\B}$, in other words, it
is the same condition as for a universal deformation with the
exception that we drop the uniqueness requirement for $f$.  If $f$
is unique whenever $\B$ is an infinitesimal algebra, then $d_{\A}$
is called a \emph{miniversal deformation} of $d$.  In \cite{ff2}, it
was shown that miniversal deformations exist whenever
$\dim(H^2(d))<\infty$. To see how to construct the miniversal
deformation, we proceed as follows.

Let $d^1=d+\delta$, where $\delta=\delta_it^i$, for some prebasis
$\{\delta_i\}$ of $H^2(d)$, the $t^i$ are parameters, and let
$\m=(t^i)$ be the ideal in $\C[[t^i]]$ generated by the $t^i$.
Denote $B=D(L_1)$ and $P_i$ be a preimage of $D(L_i)$. Note that
\begin{equation*}
[d^1,d^1]=[\delta,\delta]\in H^2\tns\m^2.
\end{equation*}
Thus
\begin{equation*}
[d^1,d^1]=\alpha_ir_2^i+\beta_2,
\end{equation*}
where $\{\alpha_i\}$ is a prebasis of $H^3$, $r_2^i\in\m^2$ and
$\beta_2\in B\tns\m^2$. Then $\beta_2=-2D(\gamma_2)$ for some
$\gamma_2\in P_2\tns\m^2$, and we define $d^2=d^1+\gamma_2$. Let
$R_2=(r_2^i)$ be the ideal generated by the $r_2^i$. Then
\begin{align*}
[d^2,d^2]&=[d^1,d^1]+2[d^1,\gamma_2]+[\gamma_2,\gamma_2]\\
&=\alpha_ir^i_2+\beta_2+2[d,\gamma_2]+2[d^1-d,\gamma_2]+[\gamma_2,\gamma_2]\\
&=\alpha_ir^i_2+2[d^1-d,\gamma_2]+[\gamma_2,\gamma_2]\\
&\in H^3\tns R_2 +L_3\tns\m^3+ L_3\tns\m^4.
\end{align*}
Moreover,
\begin{align*}
D[d^2,d^2]&=-2[d^1-d,D\gamma_2]+2[D\gamma_2,\gamma_2]\\
&=[d^1-d,\beta_2]+2[D\gamma_2,\gamma_2]\\
&=[d^1,\beta_2]+2[D\gamma_2,\gamma_2]\\
&=[d^1,[d^1,d^1]]-[d^1,\alpha_ir^i_2]+2[D\gamma_2,\gamma_2]\\
&=-[d^1-d,\alpha_ir^i_2]+2[D\gamma_2,\gamma_2]\\
&\in H^3\tns\m R_2+ L_3\tns\m^4\subseteq L_3\tns(\m R_2+\m^4)
\end{align*}
Thus
\begin{equation*}
[d^2,d^2]=\alpha_ir^i_3+\beta_3+\tau_3,
\end{equation*}
where $r^i_3-r^i_2\in\m^3$, $\beta_3\in B\tns\m^3$ and $\tau_3\in
P_3\tns(\m R^2+\m^4)$. If $R_3$ is the ideal generated by $r_3^i$,
then $\m R^2\subseteq R^3+\m^4$, so $\tau_3\in P_3\tns(R^3+\m^4)$.
Now suppose inductively that we have been able to construct $d^n$
satisfying
\begin{equation*}
[d^n,d^n]=\alpha_ir^i_{n+1}+\beta_{n+1}+\tau_{n+1},
\end{equation*}
where $r_{n+1}^i-r^i_n\in\m^{n+1}$, $\beta_{n+1}\in B\tns\m^{n+1}$,
and $\tau_{n+1}\in P\tns(R_{n+1}+\m^{n+2})$, where
$R_{n+1}=(r_{n+1}^i)$. Then $2\beta_{n+1}=-D\gamma_{n+1}$, for some
$\gamma_{n+1}\in P\tns\m^{n+1}$, and if we define
$d^{n+1}=d^n+\gamma_{n+1}$, we have
\begin{align*}
[d^{n+1},d^{n+1}]=&[d^n,d^n]+2[d^n,\gamma_{n+1}]+[\gamma_{n+1},\gamma_{n+1}]\\
=&\alpha_ir^i_{n+1}+\tau_{n+1}+2[d^n-d,\gamma_{n+1}]+[\gamma_{n+1},\gamma_{n+1}]\\
\in& H^3\tns R_{n+1}+P\tns(R_{n+1}+\m^{n+2})+ L_3\tns m^{n+2}\\
\end{align*}
Note that $L_3=H^3\oplus B\oplus P$, and thus the coboundary part in
the above lies strictly in $B\tns\m^{n+2}$. Now
\begin{align*}
D[d^{n+1},d^{n+1}]=&
[d,\tau_{n+1}]+2[D(d^n),\gamma_{n+1}]-2[d^n-d,D\gamma_{n+1}]\\&+2[D\gamma_{n+1},\gamma_{n+1}]\\
=&[d^n,\tau_{n+1}]+[d-d^n,\tau_{n+1}]+2[D(d^n),\gamma_{n+1}]+[d^n,\beta_{n+1}]
\\&+2[D\gamma_{n+1},\gamma_{n+1}]\\
=&-[d^n,\alpha_ir^i_{n+1}]+[d-d^n,\tau_{n+1}]+2[D(d^n),\gamma_{n+1}]
\\&+2[D\gamma_{n+1},\gamma_{n+1}]\\
\in& H^3\tns\m R_{n+1}+ L\tns(\m R_{n+1}+\m^{n+3})+ L\tns\m^{n+3}\\
\subseteq& L_3\tns(\m R_{n+1}+\m^{n+3}).
\end{align*}
Consequently,  we can express
\begin{equation*}
[d^{n+1},d^{n+1}]=\alpha_ir^i_{n+2}+\beta_{n+2}+\tau_{n+2},
\end{equation*}
where $r^i_{n+2}-r^i_{n+1}\in m^{n+2}$, $\beta_{n+2}\in
B\tns\m^{n+2}$ and $\tau_{n+2}\in P_3\tns(\m R_{n+1}+\m^{n+3})$.
Thus we can continue the construction indefinitely, and we obtain
finally a deformation $\dinf=d+\sum_{i=2}^\infty\gamma_i$, which
satisfies $[\dinf,\dinf]=\alpha_ir^i_\infty$, where
$r^i_\infty=\displaystyle\lim_{n\ra\infty}r^i_n$ give the relations
on the base of the miniversal deformation.

The process described above is not very efficient in constructing
the miniversal deformation, because it potentially requires an
infinite number of steps, although, in practice, it often terminates
after a finite number of steps. In some recent papers
\cite{fp3,fp4,fp6,fp7,fp8,bfp1}, the authors have shown how to
construct a miniversal deformation $\dinfty$  as follows.
\begin{equation*}
\dinfty=d+\delta_i t^i +\gamma_i x^i,
\end{equation*}
where $\{\delta_i\}$ is a prebasis of $H^2(d)$, and $\{\gamma_i\}$
is a prebasis of the 3-coboundaries $D(L_2)$, $t^i$ are the
parameters which appear in the base $\A$, and $x^i$ are formal power
series in the parameters $t^i$ which are found as follows. To
determine the coefficients $x^i$, we compute
\begin{equation*}
[\dinfty,\dinfty]=\alpha_i r^i +\beta_i s^i +\tau_i y^i,
\end{equation*}
where $\{\alpha_i\}$ is a prebasis of $H^3(d)$, $\{\beta_i\}$ is a
basis of the 3-coboundaries $D(L_2)$, and $\{\tau_i\}$ is a prebasis
of the 4-coboundaries $D(L_3)$. Thus, taken together, the
$\alpha_i,\beta_i,\tau_i$ give a basis of $L_3$. By the construction
of the miniversal deformation given in \cite{fp1}, it follows that
$s^i=0$ for all $i$, $r^i$ are formal power series in the parameters
$t^i$, and $y^i=0\mod{(r^1,\dots)}$.  Moreover, the equations
$s^i=0$ can be solved to obtain the expressions for $x^i$ as formal
power series in the parameters $t^i$. Actually, the $x^i$ can always
be expressed as rational functions of the parameters.

The construction above can be implemented on a computer, and we have
constructed Maple worksheets that carry out this implementation for
an arbitrary Lie algebra of any dimension. Using these programs we
have constructed miniversal deformations of all three dimensional
Lie algebras, which we give below.  We will also give prebases for
the cohomology for each of these examples.

\section{Calculation of the Miniversal Deformations}

\begin{table}[ht]
\begin{center}
\begin{tabular}{llrrrrr}
Type&Codiff&$H^1$&$H^2$&$H^3$\\
\hline\\
$d_1=\mathfrak{n}_3$&$\psi_1^{23}$&4&5&2\\
$d_2=\mathfrak{r}_{3,1}(\C)$&$\psi_1^{13}+\psi_2^{23}$&3&3&0\\
$d(1:1)=\mathfrak{r}_{3}(\C)$&$\psi_1^{13}+\psi_1^{23}+\psi_2^{23}$&1&1&0\\
$d(\lambda:\mu)=
\mathfrak{r}_{3,\mu/\lambda}(\C)$&$\psi_1^{13}\lambda+\psi_1^{23}+\psi_2^{23}\mu$&1&1&0\\
$d(1:0)=\mathfrak{r}_2(\C)\oplus\C$&$\psi_1^{13}+\psi_1^{23}$&2&1&0\\
$d(1:-1)=\mathfrak{r}_{3,-1}(\C)$&$\psi_1^{13}+\psi_1^{23}-\psi_2^{23}$&1&2&1\\
$d_3=\mathfrak{sl}_2(\C)$&$\psi_3^{12}+\psi_2^{13}+\psi_1^{23}$&0&0&0\\\\
\hline
\end{tabular}
\end{center}
\caption{\protect Cohomology of Three Dimensional
Algebras}\label{Table}
\end{table}
Table \ref{Table} gives the cohomology for each of the types of
3-dimensional Lie algebras. Note that for Lie algebras of type
$d(\lambda:\mu)$ there are two special values of the parameters. For
$d(1:0)$ the only variation is in $H^1$, which plays no role in the
miniversal deformation. For $d(1:-1)$ $H^2$ and $H^3$ are not the
same as for generic elements of the family $d(\lambda:\mu)$, and
this difference plays an important role in understanding the moduli
space.
\subsection{The codifferential $d_3$}
From Table \ref{Table} we see that $H^2(d_3)=0$, so there are no
deformations of this codifferential. This is not surprising, as this
codifferential corresponds to the simple Lie algebra
$\mathfrak{sl_2}(\C)$. As we shall see, some of the other
codifferentials can deform into $d_3$.
\subsection{The codifferential $d(\lambda:\mu),\mu \ne-\lambda$}
For generic values of $(\lambda:\mu)$ we have
\begin{align*}
H^1=&\langle \phd11+\phd22 \rangle\\
H^2=&\langle \psa{13}2 \rangle.
\end{align*}
The universal infinitesimal deformation is
\begin{equation*}
d^1=\psa{13}1\lambda+\psa{23}1+\psa{23}2\mu+\psa{13}2t.
\end{equation*}
Because $[d^1,d^1]=0$, the universal infinitesimal deformation is
miniversal.  Since $H^3=0$ there can be no relations on the base,
which is therefore $\A=\C[[t]]$.  The matrix of the miniversal
deformation is
\begin{equation*}
A=\left[ \begin {array}{ccc}
0&\lambda&1\\\noalign{\medskip}0&t_{{1}}&\mu\\\noalign{\medskip}0&0&0\end
{array} \right]
\end{equation*}
This matrix is just the matrix of $d(\alpha:\beta)$ where
\begin{equation*}
\{\alpha,\beta\}=\frac{\lambda+\mu\pm\sqrt{(\lambda-\mu)^2+4t}}2.
\end{equation*}
Thus deformations of $d(\lambda:\mu)$ simply move along the family.
We see that there are no jump deformations. The codifferential
$d(1:1)$ has no special behavior in this context so it is more
natural to include it in the family than $d_2$, which as we will see
later has a more complicated deformation picture.
\subsection{The codifferential $d(1:0)$}
The only special thing about $d(1:0)$ is that $H^1$ is
2-dimensional.  We have
\begin{equation*}
H^1=\langle \phd11+\phd22 \rangle.
\end{equation*}
The significance of this is that there are more outer derivations of
$d(1:0)$, which plays a role in considering extensions of a Lie
algebra by the algebra corresponding to $d(1:0)$.  The deformation
picture is generic so deformations of $d(1:0)$ simply move along the
family.  This justifies our unconventional inclusion of this element
in the one parameter family of Lie algebras.
\subsection{The codifferential $d(1:-1)$}
We have
\begin{align*}
H^1=&\langle 2\phd11+\phd21 \rangle\\
H^2=&\langle \psa{13}1,\psa{12}3 \rangle\\
H^3=&\langle\phd{123}3\rangle
\end{align*}
The universal infinitesimal deformation is
\begin{equation*}
d^1=\psa{13}1(1+t^1)+\psa{23}1-\psa{23}2+\psa{12}3t^2.
\end{equation*}
This deformation also coincides with the miniversal deformation
$\dinfty$, but in this case, we do have one relation on the base
\begin{equation*}
t^1t^2=0,
\end{equation*}
so the base of the miniversal deformation is
$\A=\C[[t^1,t^2]]/(t^1t^2)$. This relation follows from the bracket
calculation
\begin{equation*}
[\dinfty,\dinfty]=-2\phd{123}3t^1t^2.
\end{equation*}
In order to study the deformations, we have to take into account the
relation.  Thus, in the matrix
\begin{equation*}
A=\left[
\begin {array}{ccc} 0&1+t^1&1\\\noalign{\medskip}0&0&-1\\\noalign{\medskip}t^2&0&0
\end {array}
\right]
\end{equation*}
of the miniversal deformation,  we must either have $t^1=0$, or
$t^2=0$.  This means that although the tangent space $H^2$ is 2
dimensional, the actual deformations only occur along two curves.

Along the curve $t^1=0$, $\dinf\sim d_3$.  In fact, if we let $g$ be
the automorphism of $\C^3$ whose matrix $G$ is given by
\begin{equation*}
G=\left[ \begin {array}{ccc} 0&-1/2\,{\frac
{t^2-1}{t^2}}&-1/2\,{\frac {1+t^2}{t^2}}\\\noalign{\medskip}0&1&1
\\\noalign{\medskip}1&0&0\end {array} \right],
\end{equation*}
then  $d_3=g^*(\dinfty)$.  Thus we have a jump deformation from
$d(1:-1)$ to $d_3$ along this curve.

Along the curve $t^2=0$,  we have $\dinfty=d(1+t^1,-1)$, which gives
a deformation along the family.

Thus the picture for $d(1:-1)$ is more interesting.
\subsection{The codifferential $d_2$}
We have
\begin{align*}
H^1=&\langle \phd12,\phd21,\phd22 \rangle\\
H^2=&\langle \psa{13}1,\psa{13}2,\psa{23}1 \rangle.
\end{align*}
The universal infinitesimal deformation is
\begin{equation*}
d^1=\psa{13}1(1+t^1)+\psa{23}2+\psa{13}2t^2+\psa{23}1t^3.
\end{equation*}
Because $[d^1,d^1]=0$, the universal infinitesimal deformation is
miniversal, and there are no relations on the base
$\A=\C[[t^1,t^2,t^3]]$. Along the curve $t^2=t^1=0$, we have
$g^*(\dinfty)=d(1:1)$ where $g$ is determined by the matrix
\begin{equation*}
G= \left[ \begin {array}{ccc}
t_{{3}}&0&0\\\noalign{\medskip}0&1&0\\\noalign{\medskip}0&0&1
\end {array} \right].
\end{equation*}
This means that there is a jump deformation from $d_2$ to $d(1:1)$.
Otherwise,  we have that $\dinfty\sim d(\alpha:\beta)$ where
\begin{equation*}
\{\alpha,\beta\}=\frac{2+t^1\pm\sqrt{(t^1)^2-4t^2t^3}}2,
\end{equation*}
as long as $t^1\ne-2$.  Thus we see that locally (for small values
of the parameters $t^i$) $d_2$ deforms into elements of the family
$d(\lambda:\mu)$ which are ``near'' to $d(1:1)$. Thus $d_2$ deforms
as if it were the element $d(1:1)$, but is distinguishable from it
in terms of deformation behaviour because it has a jump deformation,
and in addition, it has a larger parameter space of deformations.
\subsection{The codifferential $d_1$}
We have
\begin{align*}
H^1=&\langle\phd23,\phd32,\phd11+\phd22,\phd22-\phd33\rangle\\
H^2=&\langle\psa{12}1,\psa{12}3,\psa{13}1,\psa{13}2,\psa{13}3-\psa{12}2 \rangle\\
H^3=&\langle\phd{123}2,\phd{123}3\rangle
\end{align*}
The universal infinitesimal deformation is
\begin{equation*}
d^1=\psa{23}1+\psa{12}1t^1+\psa{12}3t^2+
\psa{13}1t^3+\psa{13}2t^4+(\psa{13}3-\psa{12}2)t^5.
\end{equation*}
In this case, the miniversal deformation $\dinfty$ is equal to the
infinitesimal deformation, and we compute
\begin{equation*}
[\dinfty,\dinfty]=2\phd{123}2(t^1t^4+t^3t^5)+2\phd{123}3(t^1t^5-t^2t^3),
\end{equation*}
which means that there are two relations on the base $\A$ of the
versal deformation. We have
$\A=\C[[t^1,t^2,t^3,t^4,t^5]]/(t^1t^4+t^3t^5,t^1t^5-t^2t^3)$. The
matrix of the versal deformation is
\begin{equation*}
A=\left[ \begin {array}{ccc}
t_{{1}}&t_{{3}}&1\\\noalign{\medskip}-t_{{5}}&t_{{4}}&0\\\noalign{\medskip}t_{{2}}&t_{{5}}&0\end
{array}
 \right].
\end{equation*}
However, due to the relations, not every such $A$ is actually the
matrix of a 3 dimensional Lie algebra.  We have to solve the
relations and then consider the resulting matrices, in order to
classify the deformations. It is easy to obtain the solutions using
Maple.  We have three solutions, given by
\begin{enumerate}
\item $t^1=t^3=0$,
\item $t^3=t^4=t^5=0$,
\item $t^5=-\frac{t^1t^4}{t^3}$, $t^2=-\frac{t^1(t^4)^2}{(t^3)^2}$.
\end{enumerate}
Thus the solutions give rise to a 2-dimensional and 2 3-dimensional
pieces.

For the first solution, along the surface $t^2t^4+(t^5)^2=0$, we
have $\dinfty\sim d(1:-1)$, and otherwise $\dinfty\sim d_3$, which
means that there is a three parameter family of jump deformations to
$d_3$ and a two parameter family of jump deformations to $d(1:-1)$.

For the second solution, along the curve $t^1=\alpha+\beta$,
$t^2=\alpha\beta$, we have $\dinfty\sim d(\alpha:\beta)$, which
means that there is a jump deformation from $\dinfty$ to
$d(\alpha:\beta)$ for every value of $(\alpha:\beta)$.

For the third solution, along the surface $t^3=\alpha+\beta$,
$t^4=-\alpha\beta$, we have $\dinfty\sim d(\alpha:\beta)$, which
again gives a family of jump deformations.

As a consequence, we see that there are jump deformations from $d_1$
to every codifferential in the moduli space except $d_2$. We also
see that although the tangent space is 5 dimensional at the point
$d_1$, deformations of this codifferential actually live along lower
dimensional varieties, which illustrates a common situation in a
moduli space of Lie algebras.
\section{Conclusions}
The moduli space of 3-dimensional Lie algebras has a natural
stratification by orbifolds, three of which are just points, with
the remaining piece given by $\P^1/\Sigma_2$. The maps between these
pieces are given by jump deformations. In the case of $d_1$, there
are jump deformations to all the points in the $\P^1/\Sigma_2$
stratum.  The orbifold points play a special role in the moduli
space, either because they have extra deformations, or because there
are extra deformations to them. We illustrate the moduli space of 3
dimensional Lie algebras by the following picture.
\begin{figure}

\includegraphics[width=15cm]{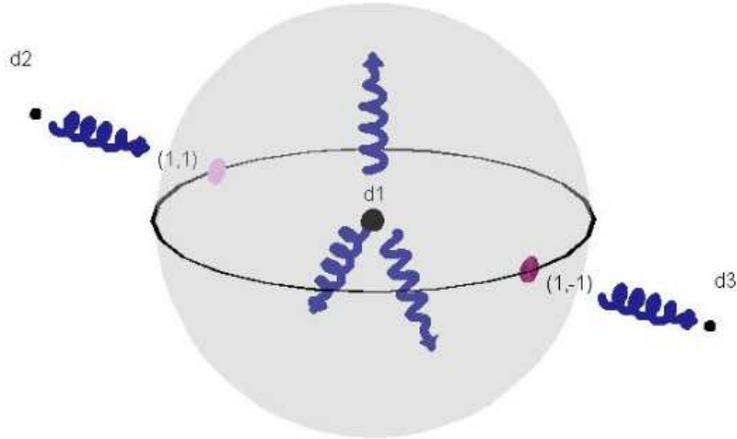}

\caption{The Moduli Space of 3 dimensional Lie Algebras}
\label{modpic}
\end{figure}

In \cite{fp8}, the moduli space of 4-dimensional Lie algebras was
studied. Although that picture is more complex, the main features,
such as the stratification by orbifolds, with jump deformations
connecting the strata, already occur in the three dimensional
picture.

\bibliographystyle{amsplain}
%\bibliography{global}
\providecommand{\bysame}{\leavevmode\hbox
to3em{\hrulefill}\thinspace}
\providecommand{\MR}{\relax\ifhmode\unskip\space\fi MR }
\providecommand{\MRhref}[2]{%
  \href{http://www.ams.org/mathscinet-getitem?mr=#1}{#2}
} \providecommand{\href}[2]{#2}

\end{document}